\documentclass{amsart}

\newtheorem*{theorem}{Theorem}
\newtheorem{cor}{Corollary}
\theoremstyle{remark}
\newtheorem*{remark}{Remark}

\DeclareMathOperator{\Ran}{Ran}
\DeclareMathOperator{\clos}{clos}
\DeclareMathOperator{\Ker}{Ker}
\DeclareMathOperator{\Int}{Int}
\DeclareMathOperator{\Ext}{Ext}
\DeclareMathOperator{\res}{res}

\begin{document}

\title{Self-adjointness of Cauchy singular integral operator}

\author{Alexey Tikhonov}

\address{Mathematical department\\
Taurida National University\\
Yaltinskaya str., 4\\
Simferopol 95007 Crimea\\
Ukraine }

\email{tikhonov@club.cris.net}


\date{}


\begin{abstract}
We extend Krupnik's criterion of self-adjointness  of the Cauchy singular integral operator to the
case of finitely connected domains. The main aim of the paper is to present a new approach for
proof of the criterion.
\end{abstract}

\maketitle

Let $G_+$ be a finitely connected domain bounded by the rectifiable curve $C=\partial G_+$,
$\,G_-=\mathbb{C}\setminus \clos\,G_+\,$ and $\infty\in G_-$. Suppose also that $w(z),\,z\in C$ is
a nonnegative weight such that $\,w(z)\not\equiv 0\,$ on each connected component of the curve $C$.
For any $\,f\in L^2(C,|dz|)\,$, we denote by $\,f_{\pm}(z)\,,\;z\in C\,$ the angular boundary
values of the Cauchy transform
\[
\mathcal{K}(f,\lambda):=\frac{1}{2\pi i}\int\limits_C\,\frac{f(z)}{z-\lambda}\,dz\,, \quad
\lambda\notin C\,
\]
from the domains $\,G_{\pm}\,$ respectively. The well-known David's theorem~\cite{Da} says that the
mappings $\,P_{\pm} \colon f \mapsto \pm f_{\pm}\,$ are bounded linear operators in
$\,L^2(C,|dz|)\,$ if and only if the curve $\,C\,$ is a Carleson curve. Moreover, the operators
$P_{\pm}$ are bounded in $\,L^2(C,w(z)|dz|)\,$ if and only if the weight $w$ is a Mackenhoupt
weight (see, e.g.~\cite{BK}), where the vector space $L^2(C,w(z)|dz|)$\, is endowed with the inner
product
\[
(f,g)_{L^2(C,w)}=\frac{1}{2\pi}\int\limits_C f(z)\,\overline{g(z)}\,w(z)\,|dz|\,,
\]
and $|dz|$ is the arc-length measure. In the sequel, we always assume that $\,C\,$ is a Carleson
curve.

In the paper we are interested in finding necessary and sufficient conditions for self-adjointness
of the projections $\,P_{\pm}\,$ and therefore (note that $\,P_++P_-=I\,$) for self-adjointness of
the corresponding Cauchy singular integral operator
\[
\mathcal{K}_S(f,\lambda):=\lim_{\varepsilon\to 0}\frac{1}{2\pi
i}\int\limits_{C(\lambda,\varepsilon)}\,\frac{f(z)}{z-\lambda}\,dz=\frac{1}{2}\,((P_+f)(\lambda)-(P_-f)(\lambda))\,,
\]
where $\,\lambda\in C\,$ and $\,C(\lambda,\varepsilon)=\{z\in C\colon |z-\lambda|<\varepsilon\}\,$.
For simple connected domains (i.e., when C is a simple closed curve) N.Krupnik~\cite{Kr} has
established the criterion: \textit{the bounded operator $\mathcal{K}_S$ is self-adjoint if and only
if $\,C$ is a circle and $w(z)\equiv const$} (see also~\cite{GK} for previous results
and~\cite{KMF} for applications). Our extension of Krupnik's criterion to the case of multiply
connected domains is a consequence of the following theorem, which is a slightly more strong
assertion than Krupnik's result.

\begin{theorem} Let $\;w\in L^1(C,|dz|)\,$.
The following conditions are equivalent:\\[3pt]
1)\,$\;\forall\,\lambda\in G_+\,\forall\,\mu\in G_-
\;\;(\frac{1}{z-\lambda},\frac{1}{z-\mu})_{L^2(C,w)}=0$;\\[3pt]
2)\, the curve $C$ is a circle and $\,w(z)\equiv const\,$.
\end{theorem}

\begin{proof} The implication $2)\,\Longrightarrow\,1)$ is obvious.\\[3pt]
$1)\,\Longrightarrow\,2)$\,.\, Let $\;w(z)\,|dz|=k(z)\,dz\,$. Obviously, $\;k\in L^1(C,|dz|)$\,.
Then $\,\forall\,\lambda\in G_+\,$
\[
 0=-\lim_{\mu\to\infty}\int\limits_C \frac{1}{z-\lambda}\cdot
\frac{\overline{\mu}}{\overline{z}-\overline{\mu}}\,k(z)\,dz=\int\limits_C
\frac{k(z)}{z-\lambda}\,dz\,.
\]
By Smirnov's theorem~\cite{Go}, we get $\,k(z)\in E^1(G_-)\,$. For the same reason, we have
$\;\forall\,\lambda\in G_+\,\forall\,\mu\in G_-\,\forall\,n\in\mathbb{N}\,$
\[
0=\int\limits_C \frac{1}{z-\lambda}\cdot
\frac{k(z)}{(\overline{z}-\overline{\mu})^{n}}\,dz\quad\Longrightarrow\quad
\frac{k(z)}{(\overline{z}-\overline{\mu})^{n}}\in E^1(G_-)\,.
\]
Let $f$ be defined by $f(z)=\overline{z}\,,\;z\in C$. Evidently, we have
\[
f(z)=\overline{\mu}+\frac{k(z)}{g(z)}\,,\quad z\in C\,,\quad\mbox{where}\quad
g(z)=\frac{k(z)}{\overline{z}-\overline{\mu}}\in E^1(G_-)
\]
and therefore the function $f(z)$ admits meromorphic continuation into the domain $G_-$. Since
$\;\forall\,n\in\mathbb{N}\,$ $\;\frac{k(z)}{(f(z)-\overline{\mu})^{n}}\in E^1(G_-)$, we get $f :
G_-\to\overline{\clos G_+}$ and $f(z)\in H^\infty(G_-)$.

The curve $C$ can be represented in the form $C=\cup_{k=0}^n C_k$, where $C_k$ are simple closed
contours. We have $\;G_+=\cap_{k=0}^n G_{k+}$, $G_-=\cup_{k=0}^n G_{k-}$, and $\infty\in G_{0-}$,
where $\;G_{0+}=\Int C_0,\; G_{0-}=\Ext C_0$, $\;G_{k+}=\Ext C_k,\; G_{k-}=\Int
C_k,\;k=\overline{1,n}$. Evidently, $f : G_{0-}\to\overline{\clos G_{0+}}$ and $f(z)\in
H^\infty(G_{0-})$. Besides, $f(z)=\overline{z}$ is an one-to-one correspondence between $C_0$ and
$\overline{C_0}=\{\bar{z}\colon z\in C_0\}$. Hence the function $f$ is a conformal mapping of
$G_{0-}$ onto $\overline{G_{0+}}$.

Let $z=\varphi(\zeta)$ be a conformal mapping of the unit disk $\mathbb{D}$ onto $G_{0+}$ and
$\zeta=\psi(z)$ be its inverse. Consider the function
$\psi_{\infty}(z)=\overline{\psi(\overline{f(z)})}$ and its inverse $z=\varphi_{\infty}(\zeta)$. It
is clear that $\,\varphi_{\infty}\,$ is a conformal mapping of $\mathbb{D}$ onto $G_{0-}$. Further,
we have
\[
\psi_{\infty}(z)=\overline{\psi(\overline{\overline{z}})}=\overline{\psi(z)}\,,\; z\in C_0\,;\quad
\psi_{\infty}(\varphi(\zeta))=\overline{\psi(\varphi(\zeta))}=\overline{\zeta}\,,\; |\zeta|=1
\]
and $\,\varphi(\zeta)=\varphi_{\infty}(\overline{\zeta})\,$. Therefore the function
\[
\Phi(\zeta)=\left\{%
\begin{array}{ll}
    \varphi(\zeta), & |\zeta|\le 1\,, \\
    \varphi_{\infty}(1/{\zeta}), & |\zeta|\ge 1\,
\end{array}%
\right.
\]
is a conformal mapping of the whole complex plane $\mathbb{C}$ onto itself and we obtain that
$\,\Phi(\zeta)=\frac{a\zeta+b}{c\zeta+d}\,$, and $\,C_0\,$ is a circle.

The curves $C_k,\;k=\overline{1,n}$ are circles too. This claim can be reduced to the case of
$C_0$: it follows easily from the observation that the operator
\[
\begin{array}{l}
C_{\varphi}\colon L^2(\varphi(C),w(\tilde{z})|d\tilde{z}|) \to L^2(C,w(\varphi(z))|dz|)\,,\\[5pt]
(C_{\varphi}f(\cdot))(z):=\sqrt{\varphi'(z)}f(\varphi(z))\,,\quad z\in C\,,\quad f\in
L^2(\varphi(C),w(\tilde{z})|d\tilde{z}|)
\end{array}
\]
is an unitary operator and from the straightforward computation
\[
(C_{\varphi}f(\cdot))(z)=\frac{c\lambda+d}{\sqrt{ad-bc}}\cdot\frac{1}{z-\lambda}\;,\quad
\tilde{z}=\varphi(z)=\frac{az+b}{cz+d}\;,\quad f(\tilde{z})=\frac{1}{\tilde{z}-\varphi(\lambda)}\;.
\]

Without loss of generality we can assume that $C_0$ is the unit circle. Other curves $\,C_k\,$ are
also circles (with centers $a_k$ and radii $r_k$). By  the same argument as above,
$\,\forall\;\mu\in G_-\; \forall\;\lambda\in G_+\quad \,$
\[
0=\int\limits_C\frac{1}{z-\mu}\cdot
\frac{\,k(z)}{\overline{z}-\overline{\lambda}}\,dz\,\quad\Longrightarrow\quad
h(z)=\displaystyle\frac{\,k(z)}{\overline{z}-\overline{\lambda}}\in E^1(G_+)\,.
\]
Since $\,k(z)\in E^1(G_-)\,$, the function $h$ admits the meromorphic continuation
$\,h(z)=\displaystyle\frac{\,k(z)}{\frac{1}{z}-\overline{\lambda}}\,$ into the domain $G_{0-}=\{z :
|z|>1\}$ and
$\,h(z)=\displaystyle\frac{\,k(z)}{\overline{a}_k+\frac{r_k^2}{z-a_k}-\overline{\lambda}}$ into the
domains $G_{k-}=\{z : |z-a_k|<r_k\}$. This meromorphic continuation has only simple poles at the
points $\,b_0=\frac{1}{\overline{\lambda}}\in G_{0-}\,$ and
$\,b_k=a_k+\frac{r_k^2}{\overline{\lambda}-\overline{a}_k}\in G_{k-}\,$. Therefore, we have
\[
h(z)-\frac{c_0}{z-b_0}-\frac{c_1}{z-b_1}\ldots-\frac{c_n}{z-b_n} \in E^1(G_+)\cap E^1(G_-)=\{0\}\,,
\]
where $\,c_k=c_k(\lambda)=\res(h(z),b_k),\;k=\overline{0,n}\,$. Hence,
\[
k(z)=c_{0}\frac{\overline{z}-\overline{\lambda}}{z-b_{0}}+
c_{1}\frac{\overline{z}-\overline{\lambda}}{z-b_{1}}+\ldots+
c_{n}\frac{\overline{z}-\overline{\lambda}}{z-b_{n}}\,,\quad z\in C\,.
\]
In particular,
\[
k(z)=-\frac{\overline{\lambda}}{z}\,\left(c_{0}+ c_{1}\,\frac{z-b_{0}}{z-b_{1}}+\ldots+
c_{n}\,\frac{z-b_{0}}{z-b_{n}}\right),\quad z\in C_0\,.
\]
Since coefficients $c_k$ depend analytically on $\,\overline{\lambda}\,$ and the functions
$\,b_k=b_k(\lambda)$ are not constants, we obtain $\,c_k=0,\;k=\overline{1,n}\,$. If $n\ne 0$, in
the same way, we have
\[
k(z)=\frac{\overline{a}_1-\overline{\lambda}}{z-a_1}\,\left(
c_{0}\,\frac{z-b_{1}}{z-b_{0}}+c_{1}+c_{2}\,\frac{z-b_{1}}{z-b_{2}}+\ldots+
c_{n}\,\frac{z-b_{1}}{z-b_{n}}\right),\quad z\in C_1\,
\]
and $\,c_k=0,\;k\ne 1\,$. Hence, $\,c_k=0,\;k=\overline{0,n}\,$ and $k(z)\equiv 0,\,z\in C$. This
contradicts our assumption $\,w(z)\not\equiv 0\,$. Thus, $\,n=0\,$,
$\,k(z)=\frac{C}{z}\,,\;|z|=1\,$ and therefore $\,w(z)\equiv const\,$.
\end{proof}
In the context of the Smirnov spaces $\,E^2(G_{\pm})\,$, the self-adjointness of the projections
$\,P_{\pm}\,$ is equivalent to the orthogonality $\,E^2(G_+)\,\bot\,E^2(G_-)\,$ (for the definition
of $\,E^2(G_{\pm})\,$, see~\cite{Go}\,). Recall that $\,\Ran P_{\pm}=E^2(G_{\pm})\,$ and $\,\Ker
P_{\pm}=E^2(G_{\mp})\,$.
\begin{cor}
If $\;E^2(G_+)\,\bot\,E^2(G_-)\;$ with respect to the
inner product $\,(\cdot,\cdot)_{L^2(C,|dz|)}\,$, then the curve $C$ is a circle.
\end{cor}
\noindent Recall that the operator $\,\mathcal{K}_S(\cdot,\lambda)\,$ are bounded in
$\,L^2(C,w(z)|dz|)\,$ if and only if $(C,w^{1/2})\in A_2$ (i.e., $w^{1/2}$ is a Mackenhoupt
weight). Under this condition we evidently have $\;w\in L^1(C,|dz|)\,$) and therefore we can
establish a desired extension of criterion of self-adjointness.
\begin{cor}
The bounded Cauchy singular integral operator $\,\mathcal{K}_S(\cdot,\lambda)\,$ is self-adjoint in
the Hilbert space $\,L^2(C,w(z)|dz|)\,$ if and only if the curve $C$ is a circle and $\,w(z)\equiv
const\,$.
\end{cor}
\noindent Note that this extension to the case of multiply connected domains can be obtained from
the corresponding result for simple connected domains. However, we prefer to regard it as
consequence of our main theorem mainly with the aim to present a new proof of Krupnik's result.
\begin{remark}
In~\cite{NF} B.Sz.-Nagy and C.Foia\c{s} employed the decomposition $\,L^2=H^2\oplus H^2_-\,$ to the
construction of their functional model for contractions, where the spaces
$\,H^2=E^2(\mathbb{D}),\,H^2_-=E^2(\mathbb{D}_-)\,$ are Hardy's spaces. We emphasize that along
with orthogonality there is an analyticity in both the domains $\,\mathbb{D}\,$ and
$\,\mathbb{D}_-\,$, respectively. However, if we intend to extend the Sz.-Nagy-C.Foia\c{s}
functional model to arbitrary domains, we cannot keep  simultaneously both these analyticities and
orthogonality because of nonorthogonality of the decomposition $\,L^2(C)=E^2(G_+)+E^2(G_-)\,$.

The combination ``analyticity only in $G_+$ plus orthogonality'' is a mainstream of development in
the multiply connected case (see, e.g.,~\cite{F} or ~\cite{BV}). In particular, the Riemann surface
(=double of the planar domain $G_+$) is used therein and the authors have to deal with the
``finite-rank defect'' decomposition $\,L^2(C)=H^2_+\oplus H^2_-\oplus \mathfrak{M}\,$, where
$0\ne\dim\mathfrak{M}<\infty$ and the subspace $\,H^2_-\,$ corresponds to the duplicate of $G_+$
(and therefore we lose entirely geometrical information concerning the domain $G_-$). Another
drawback of this approach is the use of uniformization technique or analytic vector bundles (with
fairly large amount of algebraic geometry).

On the other hand, the nonorthogonal theory keeping both analyticities and free of the above
drawbacks was developed recently in~\cite{T,Ya}. Note that the duality with respect to the Cauchy
pairing
\[
<f,g>_C:=\frac{1}{2\pi i}\int\limits_C f(z) \cdot \overline{g({\bar z})}\;dz\;,\quad f\in
L^2(C),\;g\in L^2({\bar C})
\]
is a substitute for orthogonality in our approach. For instance,
$\,E^2(G_{\pm})^{<\bot>}=E^2(\overline{G_{\pm}})\,$ (see also~\cite{TH}). Note also that linear
similarity (instead of unitary equivalence) is a natural kind of equivalence in our theory.
\end{remark}

\bibliographystyle{amsplain}

\end{document}